\documentclass[a4paper,12pt]{article}

\usepackage{amsmath,amssymb,amsthm,xypic,epsfig}

\def\g{\mathfrak{g}}
\def\h{\mathfrak{h}}
\def\n{\mathfrak{n}}
\def\b{\mathfrak{b}}
\def\C{\mathbb{C}}

\def\Z{\mathbb{Z}}
\def\qed{$\hfill \blacksquare$}
\def\l{\mathfrak{l}}

\newtheorem{theo}{Theorem}[section]
\newtheorem{prop}{Proposition}[section]
\newtheorem{lem}{Lemma}[section]

\numberwithin{equation}{section}

\begin{document}
\title{Traces of intertwiners for quantum groups and
Difference Equations II}
\date{}
\author{P. Etingof, O. Schiffmann, and A. Varchenko}
 \maketitle

\begin{abstract}
In this paper we study twisted
traces of products of intertwining operators 
for quantum affine algebras. They are interesting special functions,
depending on two weights $\lambda,\mu$, three scalar parameters
$q,\omega,k$, and spectral parameters $z_1,...,z_N$,
which may be regarded as q-analogs of 
conformal blocks of the Wess-Zumino-Witten model on an elliptic curve.
It is expected that in the rank 1 case they essentially coincide 
with the elliptic hypergeometric functions defined in \cite{FV2}. 
Our main result is that after a suitable renormalization 
the traces satisfy four systems 
of difference equations -- the Macdonald-Ruijsenaars equation,
the q-Knizhnik-Zamolodchikov-Bernard equation, and 
their dual versions. We also show that in the case when the twisting 
automorphism is trivial, the trace functions 
are symmetric under the permutation $\lambda\leftrightarrow\mu$, 
$k\leftrightarrow \omega$. Thus, our results generalize  
those of \cite{ES1}, dealing with the case $q=1$, 
and \cite{EV,ES2}, dealing with the finite dimensional case.
\end{abstract}

\section{Introduction}

Let us recall the results of \cite{EV,ES1,ES2}. 
In \cite{EV}, the main objects of study are 
suitably normalized
traces of products of intertwining operators for the quantum group
$U_q(\g)$ corresponding to a finite dimensional simple Lie
algebra $\g$. These traces are functions of two weights of $\g$, 
which take values in the endomorphism algebra of the zero weight
space of a tensor product of finite dimensional
$U_q(\g)$-modules.  
It is shown in \cite{EV} that these traces satisfy four systems of difference
equations -- the Macdonald-Ruijsenaars (MR), dual 
MR, quantum Knizhnik-Zamolodchikov-Bernard (qKZB), 
and dual qKZB equations. It is also
shown that the traces are symmetric with respect to the two
weights on which they depend. 

In \cite{ES2}, the results of \cite{EV}
are generalized to the ``twisted'' case, i.e. with traces 
involving an arbitrary diagram automorphism $T$ of $\g$. In this case,
the four systems of difference equations can still be derived,
but the symmetry no longer holds (unless $T=1$). 

It is an interesting problem to prove similar results 
in the case when $\g$ is replaced by the affine Lie algebra
$\widetilde\g$. For the classical case $q=1$, this was done earlier. 
Namely, the untwisted case $T=1$ goes back to Bernard \cite{Be}, 
while the twisted case was first considered in \cite{E1}
in a special case, and then in \cite{ES1} in general. 
In the quantum case, an attempt to attack this problem was made
in \cite{E2}, but only partial results were obtained.  

The goal of this paper is to give a complete solution of this
problem. More specifically, we derive the
Macdonald-Ruijsenaars and qKZB equations and their dual versions for 
traces of products of intertwining operators for a quantum
affine algebra, twisted by a diagram automorphism $T$ of this algebra,
and show that the trace functions are symmetric if $T=1$. 
In particular, this gives a representation-theortic 
proof of the result from \cite{FTV} that the dual difference equation 
to qKZB equation for $\g={\frak {sl}}_2$ is also the qKZB equation,
with the step and elliptic modulus interchanged,
 and allows us to generalize this result to other simple Lie algebras. 

It turns out that the methods of \cite{EV} and \cite{ES2} 
apply with only minor modifications, related to the fact that 
certain completions need to be taken in the affine case. 
Therefore, we will omit proofs of the statements, 
referring the reader to appropriate statements in \cite{EV,ES2},
and explaining only the required modifications.  

To simplify notation, we will work with simple Lie algebras 
of type A,D,E. The results can be generalized to other types. 

The structure of the paper is as follows. Section 2 serves 
to make the main definitions and set notations. In Section 3, 
we remind the basic facts about dynamical twists and R-matrices,
which are essential ingredients in the paper. In Section 4 we
introduce the trace functions -- the main characters of the paper. 
In Section 5 we discuss the MR and dual MR equations, in
Section 6 -- the qKZB and dual qKZB equations. The 
untwisted case and the symmetry identity are discussed in Section 7.

The content of this paper has the following connection with other literature. 
For $T=1$, the trace functions are q-analogs of conformal blocks 
on an elliptic curve, which are elliptic deformations 
of the q-conformal blocks on $\Bbb P^1$, defined by 
Frenkel and Reshetikhin \cite{FR}. If $\g={\frak{sl}}(2)$, 
it is expected that the trace function 
$F_T^{V_1,...,V_N}(z_1,...,z_n,\lambda,\omega,\mu,k)$
introduced here is (up to simple renomalizations) the 
elliptic hypergeometric function 
$u(z_1,...,z_n,\lambda,\tau,\mu,p)$ introduced in \cite{FV1,FV2,FTV}
by an explicit integral. 
This conjecture is motivated by the fact that 
both functions satisfy qKZB and dual qKZB equations
and the symmetry relations. The conjecture 
has been checked in the trigonometric limit 
(\cite{EV}, section 8), and in the classical limit $q\to 1$, 
where it turns into the integral formula for conformal blocks 
on the elliptic curve. We plan to discuss the general case 
in a subsequent paper.  

Another interesting case is $\g={\frak{sl}}_n$, 
and $T$ is the rotation of the affine Dynkin diagram (which is an n-gon)
by $r/n$ of a full circle, where $(r,n)=1$. 
This is the only case when the trace function does not depend 
on the dynamical parameters $\lambda$ and $\mu$, 
and is a function of $z_1,...,z_n$ and two ``elliptic moduli''
$\omega$ and $k$ (more precisely, $\mu$ takes a 
finite set of values). The qKZB equation is in this case 
the elliptic qKZ equation considered in \cite{JM,JMN},
which involves Belavin's elliptic R-matrix. Thus we show that the dual 
difference equation to this equation (i.e., its ``monodromy'')
is the difference equation involving special values of 
Felder's $R$-matrix at finitely many points. 
  
{\bf Acknowledgements.} P.E. is grateful to Igor Frenkel, 
who proposed to him the problem to derive difference equations for 
trace functions for quantum affine algebra in 1992 as a topic for a 
Ph.D. thesis. The work of P.E. was done in
part for the Clay Mathematics Institute, as a CMI Prize
fellow, and partially supported by the NSF grant 
DMS-9988796. The work of O.S. was done in part 
for the Clay Mathematics Institute. A.V. was supported by the NSF grant 
DMS-9801582.
P.E. and A.V. are grateful to MPIM for hospitality.
O.S. is grateful to MIT for hospitality.  

 \section{Notations}
\subsection{Simple 
Lie algebras} We will use the following notations:

$\g$: a simple complex Lie algebra of type A,D,E.  

$\g=\n_- \oplus \h \oplus \n_+$: its Cartan decomposition.

$\Gamma=\{\alpha_1,
\ldots,\alpha_r\}$: the Dynkin diagram.

$\Delta
 \subset \h^*$: the root system. 

$\Delta_+$: the set of positive roots. 

$(\,,\,)$: the  nondegenerate invariant bilinear form on $\g$,
normalized in such a way that 
the induced form on $\h^*$ (also denoted by $(\ ,\ )$)
satisfies the equation $(\alpha,\alpha)=2$ for roots $\alpha$.

$\theta \in \Delta^+$: the maximal root.

$\rho=\frac{1}{2}\sum_{\alpha\in \Delta^+} \alpha$: 
the half sum of positive roots.

${\rm h}$: the Coxeter (=dual Coxeter) number of $\g$.

$e_i,f_i,h_i$, $i=1,...,r$: the Chevalley-Serre generators of
$\g$. 
 
$\Omega_\h\in \h\otimes \h$: 
the inverse element to the form 
$(\,,\,)$ on $\h$.\\

\subsection{Affine 
Lie algebras}

Let $\widetilde{\g}=\g[t,t^{-1}]\oplus \C\partial
\oplus \C c$ be the affine algebra Lie associated
to $\g$, where $c$ is the central element and $\partial$ is the grading
element. That is, for Laurent polynomials $a(t),b(t)$ one has
$$
[a(t),b(t)]=[ab](t)+\text{Res}_0(a'(t),b(t))c,
$$
and $[\partial,a(t)]=ta'(t)$, $[\partial,c]=0$.
Let
$\widetilde{\b}_{\pm}$ be the Borel subalgebras of $\widetilde{\g}$ and
let $\widetilde{\n}_{\pm}$ be their nilpotent radicals.
The Cartan subalgebra of $\widetilde{\g}$ is $\widetilde{\h}=\h \oplus 
\C c \oplus\C \partial$. The dual of the Cartan subalgebra is given by 
$\widetilde{\h}^*=\h^* \oplus \C\Lambda_0 \oplus \C \delta$
where $\Lambda_0,\delta$ are defined by the relations $\langle \Lambda_0,
\h\rangle=\langle \Lambda_0,\partial\rangle=\langle\delta,\h\rangle=\langle
\delta,c\rangle=0$ and $\langle \Lambda_0,c \rangle=\langle \delta,\partial
\rangle=1$. 

\hbox to1em{\hfill} 
Set $\widehat{\g}=\g[t,t^{-1}]\oplus \C c 
=[\widetilde{\g},\widetilde{\g}]$. The Lie algebra $\widehat\g$ 
is a Kac-Moody Lie algebra.
Let $\widetilde{\Gamma}$ and $\widetilde{\Delta}$ be the
Dynkin diagram and root system of this algebra, respectively. 
Thus $\widetilde{\Delta}=
(\Delta+\Z \delta) \cup (\Z\setminus\{0\}) \delta$ and $\widetilde{\Gamma}=
\Gamma \cup
\{\alpha_0\}$ where $\alpha_0=\delta-\theta$ is the affine simple root.  
Let $(a_{ij})$ be the Cartan matrix of 
$\widetilde{\g}$.
Let 
$e_i,f_i,h_i$, $i=0,...,r$, be the standard Kac-Moody generators 
of $\widehat\g$. The span of $h_i$ (i.e. the space $\h\oplus \C c=
\widehat \g\cap \widetilde{\h}$) will   
be denoted by $\widehat \h$. 

Let $D\in \h\oplus \C\partial$ 
be the principal gradation element:
$[D,e_i]=e_i$ and $[D,f_i]=-f_i$ for all $i$. 

The Lie algebra $\widetilde\g$ carries a  
nondegenerate invariant form $(,)$, which restricts to the 
form on $\g$ defined above. 
Such a form is not unique: it depends on one free parameter.  We
will uniquely determine the form by the condition
that $(D,D)=0$. (We note that this form is somewhat different 
 from the commonly used form, which satisfies
$(\partial,\partial)=0$.) 
The restriction of this form to 
$\widetilde\h$ is nondegenerate, so it defines a form on the 
dual space. 

Let $\h'$ be the orthogonal complement 
to $\C c\oplus \C D$ in
$\widetilde\h$. It is easy to see that the map 
$x\to x-\frac{(\rho,x)}{{\rm h}}c$ defines an orthogonal 
isomorphism $\h\to \h'$. 

Let $\widetilde\rho=\rho+{\rm h}\Lambda_0$.
Then $({\widetilde\rho},\alpha_i)=1$ for $i=0,
\ldots, r$. 
Thus, $\widetilde\rho$ corresponds to $D$ under the
identification $\widetilde{\h}^*\to \widetilde{\h}$ 
defined by the inner product. 

Finally, for any $z \in \C^*$, we define an automorphism $D_z \in
 \mathrm{Aut}\;U(\widehat{\g})$ by $D_z(e_{i})=ze_{i}$ and
$D_z(f_{i})=z^{-1}f_{i}$ for any $i \in \widetilde{\Gamma}$,
and $D_zc=c$. In other words, we have $D_z(x)=z^Dxz^{-D}$. 

\subsection{Quantum groups}

Let $q$ be a nonzero complex number, such that $|q|<1$. 
We fix $\hbar$ such that $q=e^\hbar$. For any operator $A$, 
the expression $q^A$ will stand for $e^{\hbar A}$.  

Let $U_q(\widetilde{\g})$ be the quantized 
affine 
algebra corresponding to $\widetilde\g$. 
It is a Hopf algebra over $\C$ with generators 
$E_{i}$, $F_{i}$, $i=0,1,...,r$ 
and $q^{h}$, $h \in \widetilde{\h}$
and with relations
 $$q^{x+y}=q^xq^y, \
x,y\in \widetilde{\h},\quad q^0=1 
$$
$$
q^h E_{i} q^{-h}=q^{\alpha_i(h)}
E_{i},\quad
q^hF_{i}q^{-h}=q^{-\alpha_i(h)}F_{i}$$
$$E_{i}F_{j}-F_{j}E_{i}=\delta_{ij}
\frac{q^{h_{i}}-q^{-h_{i}}}{q-q^{-1}},$$
$$\sum_{k=0}^{1-a_{ij}} (-1)^k \bmatrix 1-a_{ij} \\ k
\endbmatrix_{q} E_{i}^{1-a_{ij}-k}E_{j}E_{i}^k=0,
\quad i\ne j,$$
and 
$$\sum_{k=0}^{1-a_{ij}} (-1)^k \bmatrix 1-a_{ij} \\ k
\endbmatrix_{q} F_{i}^{1-a_{ij}-k}F_{j}F_{i}^k=0,
\quad i\ne j.$$
 $$ \bmatrix n \\ k \endbmatrix
_q=\frac{[n]_q!}{[k]_q! [n-k]_q!}, \quad [n]_q! = [1]_q \cdot [2]_q
\cdot \dots \cdot [n]_q, \quad [n]_q= \frac {q^n - q^{-n}}{q-q^{-1}}$$

\paragraph{}The comultiplication $\Delta,$ antipode $S$ and counit
$\epsilon$ in $U_{q}(\widetilde{\g})$ are given by 
$$\Delta(E_{i}) = E_{i}\otimes
q^{h_{i}} + 1\otimes E_{i}, \quad
  \Delta(F_{i}) = F_{i}\otimes 1 + q^{-h_{i}}\otimes 
F_{i},$$
$$ \Delta(q^h) = q^h \otimes q^h,$$
 $$S(E_{i})=-E_{i}q^{-h_{i}},\quad
S(F_{i})=-q^{h_{i}}F_{i},\quad S(q^h)=q^{-h},$$
 $$\epsilon(E_{i}) =
\epsilon(F_{i}) = 0,\quad \epsilon(q^h) = 1.$$ 

The quantum group $U_q(\g)$ is the Hopf subalgebra of
$U_q(\widetilde\g)$ generated by 
$E_i,F_i$, $i\ge 1$, and $q^h, h\in \h$.  

If $\widetilde\lambda\in \widetilde\h^*$, then we will write 
$q^{\widetilde\lambda}$ for the element 
$q^{h_{\widetilde\lambda}}$, where 
$h_{\widetilde\lambda}$ is the image of ${\widetilde\lambda}$
under the identification $\widetilde\h^*\to \widetilde\h$. 

Let $U_q(\widehat{\g})$ be the subalgebra of
$U_q(\widetilde{\g})$
generated by $E_{i},\;F_{i}$, $\alpha_i \in \widetilde{\Gamma}$
and $q^h$, $h \in \h \oplus \C c$.
Let
$U_{q}(\widehat{\n}_{\pm})$ be the subalgebra generated by 
$(E_{i})_{\alpha_i\in \widetilde{\Gamma}}$ and
$(F_{i})_{\alpha_i\in \widetilde{\Gamma}}$ respectively, 
and $U_{q}(\widehat{\b}_{\pm})$ be generated by 
$U_{q}(\widehat{\n}_{\pm})$ and elements $q^h$,
$h\in \h\oplus \C c$. 

It is known
 that $U_{q}(\widetilde{\g})$ has a (generalized) 
quasitriangular structure, given by the R-matrix 
$$\mathcal{R} = \mathcal{R}_0q^{(c \otimes D + 
D \otimes c)/{\rm h} +\Omega_{\h'}},$$
 where $\mathcal{R}_0$ belongs to an appropriate 
completion of $U_q(\widehat\b_+)\otimes U_q(\widehat\b_-)$, 
and its degree zero part (under the principal gradation of the first
component of the tensor product) is $1$.
(Here $\Omega_{\h'}$ denotes the inverse to the form on $\h'$). 

The automorphism $D_z$ of $U_q(\widetilde{\g})$ is defined
as in the classical case. Given a finite dimensional representation $V$ of 
$U_q(\widehat{\g})$, and any $z\in \Bbb C^*$, define 
the shifted representation $V(z)$ to be the vector space $V$ with 
the new action of $U_q(\widehat{\g})$ defined by 
$\pi_{V(z)}(a)=\pi_V(D_z(a))$.

\subsection{Dynamical notation}
We will use the following notations throughout the paper. Let 
$V_1,\ldots, V_N$ be $U_q(\widehat{\g})$-modules. For an element
$A$ of $U_q(\widehat{\g})$ and 
$i\in \{1,\ldots,N\}$, we denote by $A^i$ the action of $A$ on the $i$th
component of the tensor product $V_1 \otimes \cdots \otimes V_N$. Now let
$\l \subset \widehat{\h}$ and suppose that $V_1, \ldots , V_N$ are 
$\l$-semisimple. Let $S(\lambda): \l^* \to U_q(\widehat{\g})$ be a function.
Then for any $i,j \in \{1,\ldots,N\}$ we denote by $S^i(\lambda+h^{(j)})$ the
endomorphism of $V_1\otimes...\otimes V_N$ defined by
$$S^i(\lambda+h^{(j)}) (v_1 \otimes \cdots \otimes v_N)=S_i(\lambda+\mu_j)
(v_1 \otimes \cdots \otimes v_N)$$
if $v_j$ is of weight $\mu_j \in \l^*$.

\section{Dynamical twists and quantum dynamical R-matrices}

In this section we construct dynamical twists 
and R-matrices arising from 
a Dynkin diagram authomorphism (\cite{JKOS,ESS,ES2}). 

\paragraph{}Let $T \in \mathrm{Aut}\;(\widetilde{\Gamma})$ be an automorphism
of order $N$. Let us
denote by $B_T:U_q(\widetilde{\g}) \to U_q(\widetilde{\g})$ the automorphism
defined by
$$B_T(E_i)=E_{T(i)},\quad B_T(F_{i})=F_{T(i)}, \quad
B_T(h_{i})=h_{T(i)},$$ 
and $B_T(D)=D$. For brevity, we will denote $B_T$ 
simply by $B$, assuming that $T$ has been fixed. 

It is easy to see that $B$ is an orthogonal
automorphism of $\widetilde\h$, which
preserves $\h'$,
$c$ and $D$. In particular, it preserves the
principal grading, i.e $D_zB=BD_z$ for any $z \in \C^*$.

\paragraph{}Let $\widetilde{\l} =\big(\sum_{\alpha \in 
\widetilde{\Gamma}} 
\C(\alpha-T\alpha)\big)^\perp =\widetilde\h^B
\subset \widetilde{\h}$.  
Set $\l=\widetilde{\l} 
\cap \h'$. Since $B$ is orthogonal, the form 
is nondegenerate on $\widetilde\l$ and on $\l$. We identify
$\l$ (resp. $\widetilde{\l}$) with their duals using $(\,,\,)$.
Let $\h_0 \subset \h'$
be the orthogonal complement of 
$\widetilde{\l}$ in $\widetilde{\h}$, and let $\Omega_{\h_0}
\in \h_0 \otimes \h_0$ and $\Omega_{\widetilde{\l}} \in \widetilde{\l} \otimes
\widetilde{\l}$ be the inverse
elements to the restriction of $(\,,\,)$ to $\h_0$ and $\widetilde{\l}$ 
respectively. 

Define the Cayley transform $C_T:\h_0\to \h_0$ of $T$ by 
$C_T=\frac{B+1}{B-1}$. This is a skew-symmetric operator on
$\h_0$. 

Let $I_{\pm} \subset
U_q(\widetilde{\b}_{\pm})$ be
the kernel of the projection $U_q(\widetilde{\b}_{\pm}) \to 
U_q(\widetilde{\h})$ on the elements of zero degree. Set
$$\mathcal{A}_\lambda=\C[q^{\pm 2(\lambda,\alpha_1)},
\ldots,q^{\pm 2(\lambda,\alpha_r)}].$$ We define an 
automorphism $B \in \mathrm{Aut}\;(U_q(\widetilde{\g}))$ in the same
fashion as in the classical case. Then $B D_z=D_zB$ for any $z \in \C^*$.
Let $$Z=-\frac{1}{2}((1+C_T) \otimes 1){\Omega}_{\h_0}.$$

Let $q^{-2\omega}$ be a formal parameter. More precisely, we will work
over the ring of power series in the expression $q^{-2\omega}$.

Let $\widetilde{\lambda}=\lambda+\omega\widetilde{\rho}/{\rm h}$,
$\lambda\in \l^*$. (This makes sense, since 
$\widetilde{\lambda}$ will occur only in the expression 
$q^{-2(\lambda,\alpha_i)}$, which expresses via the formal parameter 
$q^{-2\omega}$). 
 
\begin{theo}[\cite{ESS}] There exists a unique $\widetilde{\l}$-invariant 
element 
${\mathcal{J}}_{T}(\widetilde\lambda)=\sum_{m\ge 0}
{\mathcal{J}}_{T,m}(\widetilde\lambda)$, such that 
$$
\mathcal{J}_{T,0}=q^Z,
\mathcal{J}_{T,m}\in ((I_{-}[-m]\otimes I_+[m])^{\widetilde\l} 
\otimes \mathcal{A}_{\lambda})[[q^{-2\omega}]]
$$
(where $I_\pm[m]=\lbrace{x\in I_\pm, [\rho,x]=mx\rbrace})$, 
satisfying the ABRR equation (\cite{ABRR})
\begin{equation}\label{E:ABRR}
\mathcal{R}^{21}\big
( \mathrm{Ad}(q^{2\widetilde\lambda} 
B) \otimes
1\big) {\mathcal{J}}_{T}(\widetilde\lambda)=
\mathcal{J}_{T}(\widetilde\lambda)q^{\Omega_{\widetilde{\l}}}.
\end{equation} 
Moreover, ${\mathcal{J}}_{T}(\widetilde\lambda)$ satisfies the following
shifted 2-cocycle condition~:
\begin{equation}\label{E:2cocy}
\mathcal{J}^{12,3}_{T}(\widetilde\lambda)
\mathcal{J}_{T}^{12}
(\widetilde
\lambda+\frac{1}{2} h^{(3)}) =\mathcal{J}_{T}^{1,23}(\widetilde\lambda)
\mathcal{J}_{T}^{23}
(\widetilde\lambda-\frac{1}{2} h^{(1)}).
\end{equation}
\end{theo}

Below we will need a slightly modified version of $\mathcal J$. Namely, 
we set
$$
\mathbb{J}_{T}(\widetilde\lambda)=\mathcal{J}_{T}(\widetilde\lambda+
\frac{1}{2}(h^{(1)}+h^{(2)})).$$

{\bf Remark.} At $q^{-2\omega}=0$, both 
$\mathcal J_T$ and $\mathbb J_T$ specialize to $\mathcal R_0^{21}$.

Now define the exchange operator 

\begin{equation}\label{E:dynR}
\mathbb{R}_{T}(\widetilde\lambda)={\mathbb{J}}_{T}
(\widetilde\lambda)^{-1}
\mathcal{R}^{21}{\mathbb{J}}_{T}^{21}(\widetilde\lambda),
\end{equation}

\begin{prop} The element $\mathbb{R}_{T}(z,\widetilde\lambda)$ satisfies the 
\textit{quantum dynamical Yang-Baxter equation}~:
\begin{equation}\label{E:QDYBE}
\begin{split}
\mathbb{R}_{T}^{12}&(\widetilde\lambda)
\mathbb{R}_{T}^{13}(\widetilde\lambda+h^{(2)})
\mathbb{R}_{T}^{23}(\widetilde\lambda)\\
&=\mathbb{R}_{T}^{23}(\widetilde\lambda+h^{(1)})
\mathbb{R}_{T}^{13}(\widetilde\lambda)
\mathbb{R}_{T}^{12}(\widetilde\lambda+h^{(3)}).
\end{split}
\end{equation}
\end{prop}
\noindent
\textit{Proof.} By definition, $\mathbb{R}_{T}(\widetilde\lambda)^{21}$ 
is the twist of 
the quantum R-matrix $\mathcal{R}$ by the noncommutative
dynamical 2-cocycle
$\mathbb{J}_{T}(\widetilde\lambda)$. This implies the statement. 

\subsection{Spectral parameter}

\paragraph{}
Now let us introduce a spectral parameter.
Consider the expression $$\mathbb{J}_{T}(z,\widetilde\lambda)=
({D}_z \otimes 1)\mathbb{J}_{T}(\widetilde\lambda).$$
It is clear that we have 
$$\mathbb{J}_{T}(z,\widetilde\lambda) \in ((U_q(\widetilde{\b}_{-}
)
\otimes U_q(\widetilde{\b}_{+}))^{\widetilde\l}[[z^{-1}]]
\otimes \mathcal{A}_{\lambda})[[q^{-2\omega}]].$$ 
Since $B D_z=D_z B$, we have 
\begin{equation}\label{z1z2}
(D_{z_1} \otimes D_{z_2}) {\mathbb{J}}_{T}
(\widetilde\lambda)=(D_{z_1/z_2}\otimes 1){\mathbb{J}}_{T}
(\widetilde\lambda)=
{\mathbb{J}}_{T}(z_1/z_2,\widetilde\lambda).
\end{equation}

\paragraph{}Now set
\begin{equation}\label{E:dynRz}
\mathbb{R}_{T}(z,\widetilde\lambda)=(D_z \otimes 1)
\mathbb{R}_{T}(\widetilde\lambda)
\end{equation} 
Note that $\mathbb{R}_{T}(z,\widetilde\lambda)$ belongs to the space
\footnote{Here and below, 
the sign 
$\widehat{\otimes}$ will denote a completed tensor product. 
The nature of completion in each situation should be clear from the context.} 

$$ (U_q(\widetilde{\g})
{\otimes} U_q(\widetilde{\g})((z^{-1}))
\widehat{\otimes} \mathcal{A}_{\lambda})[[q^{-2\omega}]].$$

\begin{prop} The element $\mathbb{R}_{T}(z,\widetilde\lambda)$ satisfies the 
\textit{quantum dynamical Yang-Baxter equation with spectral parameters}~:
\begin{equation}\label{E:QDYBEz}
\begin{split}
\mathbb{R}_{T}^{12}&(z_{1}/z_2,\widetilde\lambda)
\mathbb{R}_{T}^{13}(z_{1}/z_3,\widetilde\lambda+h^{(2)})
\mathbb{R}_{T}^{23}(z_2/z_3,\widetilde\lambda)\\
&=\mathbb{R}_{T}^{23}(z_2/z_3,\widetilde\lambda+h^{(1)})
\mathbb{R}_{T}^{13}(z_1/z_3,\widetilde\lambda)
\mathbb{R}_{T}^{12}(z_1/z_2,\widetilde\lambda+h^{(3)}).
\end{split}
\end{equation}
\end{prop}
\noindent
\textit{Proof.} 
Equation (\ref{E:QDYBEz}) is obtained by applying $D_{z_1} \otimes
D_{z_2} \otimes D_{z_3}$ to (\ref{E:QDYBE})
and using (\ref{z1z2}).\qed

Now let $V$ be a finite dimensional representation of
$U_q(\widehat\g)$. 
Since the level of such a representation is always
zero (see e.g. \cite{CP}), we can evaluate the element
$\mathcal R$ and hence 
$\mathbb{R}_{T}(z,\widetilde\lambda)$ in the tensor product 
$V\otimes V$. The obtained operator
$\mathbb{R}_{T}^{VV}(z,\lambda,\omega)$ is a finite dimensional
solution of the quantum dynamical Yang-Baxter equation with a
spectral parameter $z$ and dynamical parameter $\lambda$, which belongs to 
$(\text{End}(V\otimes V)((z^{-1}))\otimes \mathcal{A}_\lambda)
[[q^{-2\omega}]]$.

In fact, the series defining $\mathbb{R}_T^{VV}$ has a nonempty 
region of convergence. More precisely, the 
coefficients to all powers of $q^{-2\omega}$ 
are convergent series in $z^{-1}$ 
outside a common disk, and for any point of this disk the 
resulting series in $q^{-2\omega}$ (with numerical 
coefficients) is convergent if the real part of $\omega$ is small enough.

To show this,  
for any two finite-dimensional $U_q(\widehat{\g})$-modules $V$ and
$W$ we set
$$\mathcal{R}^{V,W}(z)=(\pi_V\otimes \pi_{W})((D_z\otimes 1)\mathcal{R}).$$
This is a power series in $z$. Then we have the following
result. 

\begin{theo}[\cite{KS}]\label{KS} 
The power series $\mathcal{R}^{V,W}(z)$ converges
to a holomorphic function in the neighborhood of $z=0$, which extends to
a meromorphic function on $\Bbb C$. 
\end{theo}

Using this theorem and the Birkhoff's theory of q-difference
equations (see e.g. \cite{EFK}, Chapter 10), we get from the ABRR equation 
that the series
$\mathbb{R}^{VV}_T$ 
defines 
a meromorphic function of $z,\lambda,\omega$ in the region
$z\ne 0$, $\text{Re}(\omega)<0$.  

{\bf Remark 1.} Another proof of Theorem \ref{KS} can be found in
\cite{EM}. 

{\bf Remark 2.} We note that after evaluation in $V\otimes W$, in
the case $T=1$, the defining equation \ref{E:ABRR}
for ${\mathbb J}$ becomes
the quantum KZ equation of Frenkel and Reshetikhin \cite{FR}. 
In the case $T\ne 1$, this is a twisted version of the quantum KZ
equation.

{\bf Remark 3.} If $\g={\frak{sl}_n}$, $V$ is its vector 
representation, and $T=1$, then 
the dynamical R-matrix $\Bbb R_T^{VV}(z,\lambda,\omega)$ is 
gauge equivalent to Felder's elliptic R-matrix 
(see \cite{FR,Mo}). If $T$ is the rotation of the Dynkin diagram
of $\widehat\g$ by $r/n$ of a full circle, $(r,n)=1$, then 
$R_T^{VV}$ does not depend on $\lambda$ (as $\l=0$), and 
is a usual R-matrix with a spectral parameter. 
Namely, as was shown in \cite{T}, it is the Belavin R-matrix. 
 
\section{Twisted traces of intertwiners}

\subsection{Intertwiners}

\paragraph{}If $W$ is 
a graded space, we will denote by $W[\nu]$ the weight subspace of weight $\nu$
in $W$. 

\paragraph{}For any $\mu \in \h^{'*}$ and $k \in \Bbb C$, 
let $\widetilde\mu=\mu+k\widetilde\rho/{\rm h}$.
Let $M_{\widetilde\mu}$ be the $U_q(\widetilde{\g})$-Verma module of weight
$\widetilde{\mu}$ (it has level $k$). Let
$x_{\widetilde\mu}$ be a highest weight vector
in this module.
Let $M^*_{\widetilde\mu}$ be the restricted dual module of 
$M_{\widetilde\mu}$, and
$x^*_{\widetilde\mu}\in M^*_{\widetilde\mu}[-{\widetilde\mu}]$ 
be such that $\langle x^*_{{\widetilde\mu}},x_{\widetilde\mu}\rangle=1$.

\paragraph{}Let $V$ be a finite-dimensional
$U_q(\widehat{\g})$-module. Let $\beta\in \widehat\h^*$. 

\begin{lem}\label{FRlem}\cite{FR}
Assume that $M^*_{\widetilde\mu+\beta}$ is
irreducible. Then
$$\mathrm{Hom}_{U_q(\widehat{\g})}(M_{\widetilde\mu},M_{\widetilde\mu+\beta} 
\widehat{\otimes} V)
\simeq V[\beta]$$
where the isomorphism is as follows : to $\Phi : M_{\widetilde\mu} \to
M_{\widetilde\mu+\beta} \widehat{\otimes} V$ we associate the element 
$\langle \Phi \rangle=
\langle x^*_{{\widetilde\mu+\beta}}, \Phi(x_{\widetilde\mu})\rangle$.
\end{lem}

 We will denote the operator 
corresponding to $v\in V[\beta]$ under the correspondence
of Lemma \ref{FRlem} by $\Phi_{\widetilde\mu}^v$. 
Define also
the intertwining operator $\Phi_{\widetilde\mu}^v(z):
M_{\widetilde\mu}\to M_{\widetilde\mu+\beta}\widehat{\otimes}V(z)$. 
This operator can be regarded as an infinite in both directions 
Laurent series 
in $z$ with operator coefficients
(for the coefficients, the completion of the tensor product is 
unnecessary).  

\paragraph{}Let $V_1,\ldots, V_N$ be finite-dimensional
$U_q(\widehat{\g})$-modules and let $v_1 \in V_1[\nu_1],\ldots$,\\
$v_N\in V_N[\nu_N]$. For generic values of $\mu \in \h^{'*}$ and $k$
we can consider the
composition
$$
\Phi^{v_1,\ldots v_N}_{{\widetilde\mu}}(z_1, \ldots, z_N)=
\Phi^{v_1}_{{\widetilde\mu}-\nu_2-\ldots-\nu_{N}}(z_1)\cdots
\Phi^{v_N}_{{\widetilde\mu}}(z_N) :$$
$$M_{\widetilde\mu}\to M_{{\widetilde\mu}-\nu_1-\ldots-\nu_N} 
\widehat{\otimes} V_1\widehat{\otimes} \cdots \widehat{\otimes}
V_N,
$$
It is easy to see that, for a given $\tilde\mu$,
the matrix coefficients of $\Phi^{v_1,\ldots,v_N}_{\widetilde\mu}(z_1,\ldots,z_N)$ belong
to the space $z_1^{l_1}\cdots z_N^{l_N}
\C[[\frac{z_2}{z_1},\ldots,\frac{z_{N}}{z_{N-1}}]]$ for some $l_1,\ldots,l_n
\in \Z$.
\begin{theo}[\cite{FR},\cite{EFK}]\label{T:01}
The matrix elements of the operator valued 
series $\Phi^{v_1,\ldots,v_N}_{\widetilde\mu}
(z_1,\ldots,z_N)$ converge in the the region $z_1 \gg z_2\gg \cdots \gg z_N$
and extend to meromorphic 
functions on $(\C^*)^N$.\end{theo}

\subsection{Trace functions}
\paragraph{}For any $\mu \in \h^{'*}$ we define the linear operator $B:
M_{\widetilde\mu}\to M_{B\widetilde\mu}$ by setting $B(u\cdot x_{\widetilde\mu})=B(u)\cdot x_{B\widetilde\mu}$ for any $u\in
U_q(\widetilde{\n}_{-})$.
\paragraph{}Consider the following {\it trace function}, 
which will be the main
object of study of this paper.
$$\Psi^{v_1,\ldots,v_N}_{T}(z_1,\ldots,z_N,\lambda,\omega,\mu,k)=
\mathrm{Tr}_{|M_{\widetilde\mu}}(\Phi^{v_1,\ldots,v_N}_{{B\widetilde\mu}}(z_1,\ldots,z_N)B
q^{2\widetilde\lambda}),$$
where $\widetilde\lambda=\lambda+\omega\widetilde\rho/{\rm h}$, 
and $\lambda \in {\l}^*$. Note that 
$\Psi^{v_1,\ldots,v_N}_{T}(z_1,\ldots,z_N,\lambda,\omega,\mu,k)=0$ unless
$\sum_i \nu_i=B{\mu}-{\mu}$, and that 
$\Psi$ is independent of $z_i$ if $N=1$. Moreover, given
$\nu_1,\ldots,\nu_N$, the set of all $\mu \in \h^{'*}$ satisfying this condition
is an $\l^*$-homogeneous space. By Theorem~\ref{T:01}, for (almost) every such
fixed value of $\mu$ and $k$, the series
$\Psi^{v_1,\ldots,v_N}_{T}(z_1,\ldots,z_N,\lambda,\omega,\mu,k)$ converges
in the region $z_1 \gg \cdots \gg z_N$ and extends
to a meromorphic function on $(\C^*)^N$ with values in 
$q^{2(\lambda,\mu)}\mathcal{A}_\lambda[[q^{-2\omega}]]$. 

Now we will consider  $k$
as a formal parameter. Namely, we set
$$\mathcal{A}_\mu^-=\C[q^{-2(\mu,\alpha_1)},
\ldots,q^{-2(\mu,\alpha_r)}],\qquad \mathcal{A}_{\mu,k}=
\mathcal{A}_\mu^-[[q^{-2k}]].$$
The following lemma is easy to prove.
\begin{lem} The function 
$\Psi^{v_1,\ldots,v_N}_{T}(z_1,\ldots,z_N,\lambda,\omega,\mu,k)$
has a series expansion, which lies in 
$(V_1 \otimes \cdots \otimes V_N)^\l
\otimes q^{2(\lambda,\mu)}
\mathcal{A}_{\lambda,\omega}\widehat{\otimes}\mathcal{A}_{\mu,k},$
where $\widetilde{\mu}=\mu+k{\widetilde\rho}/{\rm h}$ and $\widetilde{\lambda}=\lambda+ \omega
{\widetilde\rho}/{\rm h}$.
\end{lem}

\paragraph{}We also define the following universal trace function, which
takes values in
$(V_1 \otimes \cdots \otimes V_N)^\l \otimes (V_N^*
\otimes\cdots \otimes
V_1^*)^\l\otimes q^{2(\lambda,\mu)}\mathcal{A}_{\lambda,\omega}
\widehat{\otimes}\mathcal{A}_{\mu,k}$ :
\begin{equation*}
\begin{split}
\Psi^{V_1,\ldots,V_N}_{T}&(z_1,\ldots,z_N,\lambda,\omega,\mu,k)\\
&=\sum_{v_i \in \mathcal{B}_i}
\Psi^{v_1,\ldots,v_N}_{T}(z_N,\ldots,z_1,\lambda,\omega,\mu,k)\otimes
v_N^* \otimes\cdots \otimes v_1^*
\end{split}
\end{equation*}
where $\mathcal{B}_i$ is a homogeneous basis of $V_i$ and 
$\{v_i^*\}_{v_i \in \mathcal{B}_i}$ is the dual basis of $V_i^*$.
This function is defined for $\mu\in \h^{'*}$, but 
has support on finitely many shifted copies of $\l^*$. 

\paragraph{}It is more convenient to renormalize the function
$\Psi^{V_1,\ldots,V_N}_{T}$ as follows.
Define
$$\mathbb{J}_{T}^{1\cdots N}(\widetilde\lambda)=
\mathbb{J}_{T}^{1,2\cdots N}
(\widetilde\lambda)\cdots \mathbb{J}_{T}^{N-1,N}(\widetilde\lambda),$$
$$\mathbb{J}_{T}^{1\cdots N}(z_1,\ldots,z_N,\widetilde\lambda)=
(D_{z_1} \otimes \cdots \otimes D_{z_N})
\mathbb{J}_{T}^{1\cdots N}(\widetilde\lambda).$$
It is easy to check that
$$\mathbb{J}_{T}^{1\cdots N}(z_1,\ldots,z_N,\tilde\lambda)
\in U_q(\widetilde{\g})^{\otimes N}[[\frac{z_2}{z_1},\ldots,
\frac{z_N}{z_{N-1}}]]\widehat{\otimes} \mathcal{A}_{\lambda,\omega}.$$

Further, 
put $\mathbb{Q}_{T}(\widetilde\lambda)=m_{21}((1 \otimes S)[
(\mathcal R^{21})^{-1}\mathbb{J}_{T}
(\widetilde\lambda)])$ (here $m_{21}$ is the multiplication of components, 
i.e. $m_{21}(a\otimes b):=ba$). 
This is a formal power series lying in 
$U_q(\widetilde{\g})
\widehat{\otimes} \mathcal{A}_{\lambda,\omega}$.
Note that ${\mathbb{Q}}_{T}(\widetilde\lambda)$ is $B$-invariant and
$D_z$-invariant. If $T=1$, we will drop the subscript $T$ and 
denote the corrseponding elements $\Bbb J_T,\Bbb Q_T$, $\Bbb R_T$
simply by $\Bbb J$, 
$\mathbb R$, $\mathbb{Q}$. 

Finally, let
$$\delta_q^{T}({\lambda},\omega)=\big(\mathrm{Tr}_{|M_{-\widetilde{\rho}}}(
Bq^{2(\lambda+ \omega{\widetilde\rho}/{\rm h})})\big)^{-1}$$
be the twisted Weyl denominator.
It can be written as an infinite product of binomial factors
(see \cite{ES1}). 

\paragraph{}Define the renormalized trace function by
\begin{equation*}
\begin{split}
&F^{V_1,\ldots ,V_N}_{T}(z_1,\ldots,z_N,\lambda,\omega,\mu,k)\\
&=[{\mathbb{Q}}^{-1}(\mu+h^{(*1\cdots *N)},k)
^{(*N)}\otimes \cdots \otimes {\mathbb{Q}}^{-1}
(\mu+h^{(*1)},k)^{(*1)}]\\
&\qquad\qquad\qquad\qquad
\varphi^{V_1,\ldots,V_N}_{T}(z_1,\ldots,z_N,\lambda,\omega,-\mu-\rho
,-k-{{\rm h}})
\end{split}
\end{equation*}
where
\begin{equation*}
\begin{split}
\varphi_{T}^{V_1,\ldots, V_N}&(z_1,\ldots,z_N,\lambda,\omega,\mu,k)\\
&=\mathbb{J}_{T}^{1\cdots N}(z_1,\ldots,z_N,\lambda,\omega)^{-1}
\Psi_{T}^{V_1,\ldots, V_N}(z_1,\ldots,z_N,\lambda,\omega,\mu,k)
\delta_q^{T}(\lambda,\omega).
\end{split}
\end{equation*}
(here $\Bbb Q(\lambda,\omega):=\Bbb Q(\widetilde\lambda)$,
$\Bbb J(\lambda,\omega):=\Bbb J(\widetilde\lambda))$). 
Note that $F^{V_1,\ldots ,V_N}_{T}(z_1,\ldots,z_N,\lambda,\omega,\mu,k)$
is in the space
$$\mathcal{E}=(V_1\otimes \cdots \otimes V_N)^\l\otimes (V_N^*
\otimes \cdots \otimes V_1^*)^\l
((\frac{z_2}{z_1},\ldots,\frac{z_{N}}{z_{N-1}}))
\widehat{\otimes} q^{-2(\lambda,\mu)}
\mathcal{A}_{\lambda,\omega} \widehat{\otimes} \mathcal{A}_{\mu,k}.$$

{\bf Remark 1.} Here we regard the function $F_T$ as a formal power
series. It is expected, but, as far as we know, 
not proved in general, that this series is actually convergent 
in some region, defining an analytic function.
In the case $T=1$, it is expected
that $F_T$ has an integral representation coming from the free 
bosonic realization, which could be used for a proof 
of analyticity. Such a presentation has not been worked
out in general. However, in the case of $\g={\frak{sl}}_2$, 
such presentation is known (\cite{Ko1,Ko2}), which 
means that in this case the function $F_T$ for $T=1$ is analytic.

{\bf Remark 2.} It is important to note that the element
$\Bbb Q_T(\widetilde\lambda)$ defined here 
is somewhat different from the one used in \cite{EV,ES2}.
This is done because the formula 
$m_{21}((1\otimes S^{-1})(\Bbb J_T(\lambda)))$ used in \cite{EV,ES2}
does not make sense for affine Lie algebras 
due to appearance of infinite summations. Indeed, 
at $q^{-2\omega}=0$ we have $\Bbb J=\mathcal R^{21}_0$, so 
already the computation of the leading term of 
$\Bbb Q_T$ as in \cite{EV,ES2} in a finite dimensional 
representation $V$  
would involve evaluating the R-matrix 
of $V,V^*$ at $z=1$ and multiplying the dual to the second component 
by the first component. However, this operation contains 
infinite summations, which are manifestly illegal, since the R-matrix
$\mathcal R_{V,V^*}(z)$ has a pole at $z=1$. 
On the other hand, the definition of $\Bbb Q_T(\lambda)$ given 
here makes sense, because the expression 
$(\mathcal R^{21}_0)^{-1}\Bbb J$ is a Taylor series
in $q^{-2\omega}$, whose coefficients involve only 
 finite summations. 

Now we explain the connection between our definition 
and that of \cite{ES2} in more detail. 

In the finite dimensional case, both definitions make sense, 
and the relationship 
is given by the formula
$\Bbb Q_T(\lambda)=S(u)^{-1}\Bbb Q_T^{ES2}(\lambda)$, 
where $u$ is the 
Drinfeld element \cite{D}, $u=m_{21}((1\otimes S)(\mathcal R))$,
The trace functions $F_T^{V_1,...,V_N}$ 
can be defined using either definition:
the difference would be just a constant factor.   

In the affine case, neither of the factors $S(u)^{-1}$ and    
$\Bbb Q_T^{ES2}(\widetilde\lambda)$ makes sense 
because of infinite summations,
but the product $S(u)^{-1}\Bbb Q_T(\tilde\lambda)$ does. 
So we use this product as a 
substitute for the non-existenct $\Bbb Q_T^{ES2}(\widetilde\lambda)$. 

{\bf Remark 3.} One may also consider traces
over integrable modules rather than Verma modules. 
In the particular case $\g={\frak{sl}}_n$, 
$N=1$, $V_1=S^{kn}\Bbb C^n$ they are affine 
analogues of Macdonald's polynomials and were considered in 
\cite{EK1}. 

\section{The twisted Macdonald-Ruijsenaars equations}

\subsection{The equations}

\paragraph{}Let $W$ be an integrable 
lowest weight $U_q(\widetilde{\g})$-module
of level $\eta$ (so $\eta$ is a nonpositive integer). 
Consider the following difference operator with respect to the variable
$\lambda$ acting on functions in the space $\mathcal{E}$:
\begin{equation}\label{E:0pi1}
\begin{split}
&\mathcal{D}^{T}_W\\
&\;=\sum_{\nu\in \l^*,a\in \Bbb C} 
\mathrm{Tr}_{|W[\nu+a\delta+\eta\widetilde\rho/{\rm h}]}
\big(\mathbb{R}_{T}^{WV_1}(z_1^{-1},
\lambda +h^{(2\cdots N)},\omega)
\cdots \mathbb{R}_{T}^{WV_N}(z_N^{-1},\lambda
,\omega)\big)
q^{-2ka}\mathbb{T}_\nu
\end{split}
\end{equation}
where $\mathbb{T}_\nu f(\lambda)=f(\lambda+\nu)$,
and $\Bbb R(z,\lambda,\omega):=\Bbb R(z,\widetilde\lambda)$.
In the above, we only
consider the trace of the ``diagonal block'' of the product
$\mathbb{R}_{T}^{WV_1}
(z_1^{-1},\lambda + h^{(2\cdots N)},\omega)\cdots \mathbb{R}_{T}^{WV_N}
(z_N^{-1},\lambda,\omega)$,
i.e the part that preserves $W[\nu+a\delta+\eta\widetilde\rho/{\rm h}]$.
\begin{theo}[Macdonald-Ruijsenaars equations]\label{T:02}
\begin{equation}\label{E:MR}
\begin{split}
\mathcal{D}_W^{T} F_{T}^{V_1,\ldots, V_N}&
(z_1,\ldots z_N,\lambda,\omega+\eta,\mu,k)\\
&=\chi_W(q^{-2{\widetilde{\mu}}})
F_{T}^{V_1,\ldots, V_N}
(z_1,\ldots,z_N,\lambda,\omega,\mu,k),
\end{split}
\end{equation}
where $\chi_W(x)=\sum \mathrm{dim}\;W[\nu]x^\nu$ is the character of $W$ and
where $\mathcal{D}_W^{T}$ acts on the variable $\lambda$.\end{theo}

{\bf Remark 1.} Note that equation (\ref{E:MR}) is an infinite 
difference equation. That is, the difference operator on the left 
hand side involves an infinite 
linear combination of shifts by different weights.
Nevertheless, it makes sense as an operator on power series 
in the formal parameters $q^{-2\omega}$ and $q^{-2k}$, 
since the coefficient of any fixed power of 
$q^{-2k}$ is an honest (finite) difference operator.
Such operators were previously considered in \cite{E3}.
Abusing terminology, we will call them ``difference operators'',
without specifying that they are infinite. 

{\bf Remark 2.}
We do not know a simple explicit formula for the operator \ref{E:0pi1}
in any nontrivial case.  

\paragraph{}Let us denote by $W^B$ the twist of $W$ by $B$: as a vector space
$W=W^B$ 
and the $U_q(\widetilde{\g})$-action is given by $\pi_{W^B}(u)
=\pi_W(B^{-1}(u))$.
Now suppose that $W\simeq W^B$ as $U_q(\widetilde{\g})$-modules and let us
consider the intertwiner $W\to W^B$ which preserves the lowest weight vector 
(abusing notation, we will denote it by the same letter $B$). 
This endows $W$ with the structure of a
 module over
$\Bbb C[\langle B \rangle] \ltimes U_q(\widetilde{\g})$,
where $\langle B \rangle$ is the cyclic group generated by $B$. 

\paragraph{}Consider the following difference operator with respect to the
variable ${\widetilde\mu}$ acting on functions in
the space $\mathcal{E}$:
\begin{equation}
\begin{split}
 &\mathcal{D}_W^{\vee,T}\\
&=\sum_{\nu,a} \mathrm{Tr}_{|W[\nu+a\delta+\eta\widetilde\rho/{\rm h}]}\big(\mathbb{R}^{WV^{*}_N}(z_N^{-1},
{\mu} + h^{(*1\cdots *N-1)},k)
\cdots \mathbb{R}^{WV^*_1}(z_1^{-1},\mu,k)
B_{W}\big)q^{-2\omega a}\mathbb{T}^\vee_\nu
\end{split}
\end{equation}
 where the sign $\vee$ means that the operators act on functions of 
$\mu$ rather than $\lambda$. 
\begin{theo}[Dual Macdonald-Ruijsenaars equations]\label{T:03} 
\begin{equation}\label{E:DMR}
\begin{split}
\mathcal{D}^{\vee,T}_{W}F_T^{V_1,\ldots,V_N}&(
z_1,\ldots,z_n,\lambda,\omega,\mu,k+\eta)\\
&=\mathrm{Tr}_{|W}(q^{-2\widetilde{\lambda}}B)
F_{T}^{V_1,\ldots, V_N}(z_1,\ldots,z_N,\lambda,\omega,\mu,k).
\end{split}
\end{equation}
\end{theo}
\paragraph{Remark.}Any $B$-invariant integrable 
lowest weight $U_q(\widetilde{\g})$-module 
is a direct sum of modules $\overline{V}^-_{\widetilde{\nu}_0}:
=\bigoplus_{\widetilde{\nu} \in \langle B
 \rangle \widetilde{\nu}_0} 
V^-_{\widetilde{\nu}}$, where $V^-_{\widetilde{\nu}}$ 
is the irreducible lowest weight module 
of lowest weight $-\widetilde{\nu}$, and  
$\widetilde{\nu}_0$ is a dominant integral weight. It is easy to 
see that both sides of (\ref{E:DMR}) identically vanish when
$W=\overline{V}_{\widetilde{\nu}_0}$ and $\widetilde{\nu}_0 \not\in \l^*$ 
(i.e when $B(\widetilde{\nu}_0) \neq \widetilde{\nu}_0$). 

\subsection{The proofs} 
\paragraph{}The proof of Theorem~\ref{T:02} is an extension of the proofs of
Theorem 2.1 of \cite{ES2} and Theorem 1.1 of \cite{EV}.
We first introduce the notion of (twisted) radial part.
\begin{prop}(see \cite{E3}) Let $V$ be a finite-dimensional
$U_q(\widehat{\g})$-module. For any element $X \in U_q(\widehat{\g})$
there exists a unique difference operator (with respect to the variable
$\lambda$) $\mathcal{D}_X^T$ acting on power series in the space
$V^{{\l}} \otimes
q^{2(\lambda,\mu)}\mathcal{A}_{\lambda,\omega}$
 such that we have
$$\mathrm{Tr}_{|M_\mu}(\Phi^V_{B\mu} X q^{2\widetilde{\lambda}}B)=
\mathcal{D}_X^T \mathrm{Tr}_{|M_\mu}(\Phi^V_{B\mu} q^{2\widetilde{\lambda}}B).$$
\end{prop}

The operator $\mathcal{D}_X^T$ is called the twisted radial part of $X$.

\paragraph{}Let $W$ be any lowest-weight integrable
$U_q(\widetilde{\g})$-module. Consider the expression
$$C_W=\mathrm{Tr}_{|W}(1 \otimes \pi_W)(\mathcal{R}^{21}\mathcal{R}(1
 \otimes q^{2\widetilde{\rho}})).$$
Set $\mathcal{M}^T_W=\mathcal{D}^T_{C_W}$.
\begin{prop} We have
\begin{enumerate}
\item $\mathcal{M}_W^T \mathcal{M}_{W'}^T=\mathcal{M}_{W'}^T \mathcal{M}_W^T$
for any two integrable modules $W,W'$,
\item For any finite-dimensional $U_q(\widehat{\g})$-module $V$ we have
$\mathcal{M}^T_W \Psi^V_{T}(z,\lambda,\omega,\mu,k)=
\chi_W(q^{2(\widetilde{\mu}+\widetilde\rho)}) 
\Psi_{T}^V(\lambda,\omega,\mu,k)$ where 
$\chi_W(x)=\sum_\nu \mathrm{dim}\;W[\nu]
x^\nu$ is the character of $W$.
\end{enumerate}
\end{prop}
\noindent
\textit{Proof.} See \cite{EK}, \cite{EV}.
\paragraph{}The computation of the operator $\mathcal{M}^T_W$ is now
identical to that in \cite{ES2}, Section 3. We obtain
\begin{equation}\label{E:0*}
\begin{split}
\mathcal{M}^T_W \Psi^T_V&(\lambda,\omega,\mu,k)=
\sum_{\nu,a} \bigg\{q^{2(k+{{\rm h}})a}\frac{\delta^{T}_q(\lambda+\nu,\omega+\eta)}{\delta_q^{T}
(\lambda,\omega)}\times \\ &
\mathrm{Tr}_{|W[\nu+a\delta+\eta
\widetilde\rho/{\rm h}]}
(\mathbb{R}_{T}^{WV}(z^{-1},\lambda,\omega))\Psi^V_{T}
(z,\lambda+\nu,\omega+\eta,\mu,k)\bigg\}.
\end{split}
\end{equation}
Theorem~\ref{T:02} now follows from (\ref{E:0*}) and from the fusion identity

\begin{equation}
\mathbb{J}_{T}^{1\cdots N}(\widetilde\lambda)^{-1}
\mathbb{R}_{T}^{0,1\cdots N}(\widetilde\lambda)
\big)\mathbb{J}_{T}^{1\cdots N}(\widetilde\lambda+h^{(0)}=
\mathbb{R}_{T}^{01}
(\widetilde\lambda+h^{(2\cdots N)})\cdots(
\mathbb{R}_{T}^{0N}(\widetilde\lambda)).
\end{equation}

\paragraph{}The proof of Theorem~\ref{T:03} is completely parallel to the
proofs of Theorem 2.3 in \cite{ES2} and Theorem 1.2 in \cite{EV}.

\section{The twisted qKZB and dual qKZB equations}

\paragraph{}Let $V_1,\ldots,V_N$ be irreducible finite-dimensional
$U_q(\widehat{\g})$-modules.
For each $j \in \{1,\ldots,N\}$ we define the
following operators acting on the space $\mathcal{E}$ :
\begin{align}
D^{T}_j=&q_{*j}^{-2\mu-C_{\h'}}q_{*j,*1}^{-2\Omega_{\h'}}\cdots
q_{*j,*j-1}^{-2\Omega_{\h'}},\\
K^{T}_j=&\mathbb{R}_{T}^{j+1,j}(\frac{z_{j+1}}{z_j},\lambda
+h^{(j+2,\ldots, N)},\omega)^{-1} \cdots \mathbb{R}^{Nj}_{T}(\frac{z_N}{z_j},
\lambda,\omega)^{-1}\Gamma_j \times\notag\\
&\mathbb{R}_{T}^{j1}(\frac{z_j}{z_1},\lambda+
h^{(2\ldots, j-1)}+h^{(j+1\ldots, N)},\omega)
\times\cdots \mathbb{R}_{T}^{j,j-1}(\frac{z_j}{z_{j-1}},
\lambda+h^{(j+1 \ldots, N)},\omega)
\end{align}
where $C_{\h'}= m_{12}(\Omega_{\h'})\in U(\h')$ 
is the quadratic Casimir element for
$\h'$, and where $\Gamma_jf(\lambda)=f(\lambda+h^{(j)})$.
\begin{theo}[qKZB equations] For all $j=1,\ldots, N$, we have :
\begin{equation}\label{E:qKZB}
\begin{split}
F_{T}^{V_1,\ldots, V_N}&(z_1,\ldots,pz_j,\ldots,z_N,\lambda,\omega,\mu,k)\\
&=(D^{T}_j \otimes
K^{T}_j) F^{V_1,\ldots, V_N}_{T}(z_1,\ldots,z_N,\lambda,\omega,\mu,k),
\end{split}
\end{equation}
where $p=q^{2k}$.
\end{theo}
\paragraph{}Similarly, consider the operators
\begin{equation}\label{E:71}
\begin{split}
D^{\vee,T}_j=&q_j^{-2\lambda-C_\l}q_{j,j+1}^{-\Omega_\l}\cdots
q_{j,N}^{-\Omega_\l},\\
K^{\vee,T}_j=&\mathbb{R}^{*j-1,*j}(\frac{z_{j-1}}{z_j},
\mu+h^{(*1\cdots *j-2)},k)^{-1}\cdots 
\mathbb{R}^{*1,*j}(\frac{z_1}{z_j},
\mu,k)^{-1}\Gamma^*_{B^{-1}(j)}\times\\
&\mathbb{R}^{*j,*N}(\frac{z_j}{z_N},\mu+h^{(*j+1\cdots *N-1)}+
h^{(*1\cdots *j-1)},k)\times\cdots\times\\
&\mathbb{R}^{*j,*j+1}(\frac{z_j}{z_{j+1}},\mu+h^{(*1\cdots *j-1)},
k),
\end{split}
\end{equation}
where $C_\l=m_{12}(\Omega_\l) \in U(\l)$ and where $(\Gamma^*_{B^{-1}(j)}
f)(\mu)=f(\mu+B^{-1}(h^{(*j)}))$.
\begin{theo}[Dual qKZB equations]\label{T:04}
For all $j=1\ldots, N$, we have :
\begin{equation}\label{E:DqKZB}
\begin{split}
F^{V_1,\ldots, V_j,\ldots, V_N}_{T}(z_1,\ldots,
sz_j,\ldots,z_N,\lambda,
,\omega,\mu,k)\\
=(D^{\vee,T}_j\otimes K^{\vee,T}_j) F_{T}^{V_1,\ldots,
V_j^B,\ldots, V_N}(z_1,\ldots,z_N,\lambda,\omega,\mu,k),
\end{split}
\end{equation}
where $s=q^{2\omega}$, and we identify $V_j$ with $V_j^B$.
\end{theo}

\paragraph{}The proof of Theorem 5.1 is very similar to the proof of Theorem
2.2 of \cite{ES2} and Theorem 1.2 of \cite{EV}
(with modifications described in Remark 2 in Section 4).
The proof of the dual qKZB equation runs parallel to \cite{ES2},
Section 6.

\section{The symmetry identity}

\paragraph{}In this section we assume that $T=1$. Thus $\l=\h'$ and the
renormalized trace function
$F^{V_1,\ldots ,V_N}(z_1,\ldots,z_N,\lambda,\omega,\mu,k)$
takes values in the space
$$(V_1\otimes \cdots \otimes V_N)[0]\otimes (V_N^*
\otimes \cdots \otimes V_1^*)[0]
((\frac{z_2}{z_1},\ldots,\frac{z_{N}}{z_{N-1}}))\widehat{\otimes} 
q^{-2(\lambda,\mu)}
\mathcal{A}_{\lambda,\omega} \widehat\otimes \mathcal{A}_{\mu,k}.$$
As in the case of
finite-dimensional simple Lie algebra $\g$ (see \cite{EV}), we have the
following symmetry identity :
\begin{theo}
$$F^{V_1,\ldots,V_N}(z_1,\ldots,z_N,\lambda,\omega,\mu,k)=
F_*^{V_N^*,\ldots,V_1^*}(z_N,\ldots,z_1,\mu,k,\lambda,
 \omega),$$
where $F_*$ is the result of interchanging the two factors
$(V_1 \otimes \cdots \otimes V_N)[0]$ and $(V_N^* \otimes \cdots \otimes
V_1^*)[0]$.
\end{theo}
This result follows from Theorem 4.1 and Theorem 4.2 and the arguments in
\cite{EV}, Section 5.

\small{

\noindent Pavel Etingof, MIT Mathematics Dept., 77 Massachusetts Ave.,
 Cambridge 02139 MA., USA\\
\texttt{etingof@math.mit.edu}\\
\noindent Olivier Schiffmann, Yale Mathematics Dept., 10 Hillhouse Ave.,
 New Haven, 06510 CT., USA\\
\texttt{schiffmann@math.yale.edu}\\
\noindent Alexander Varchenko, Department of Mathematics, University of North 
Carolina, Chapel Hill, NC 27599, USA\\
\texttt{anv@email.unc.edu}
 \end{document}